\documentclass[12pt]{amsart}
\usepackage{times}
\usepackage{amsfonts, amsmath, amssymb, amsthm, mathtools, dsfont}
\usepackage[margin=3cm]{geometry}
\usepackage[colorlinks=true,citecolor=blue]{hyperref}
\usepackage{setspace, enumerate}
\usepackage{comment}
\usepackage[title]{appendix}

\theoremstyle{plain}
\newtheorem{theorem}{Theorem}[section]

\newtheorem{lemma}[theorem]{Lemma}

\newtheorem{claim}[theorem]{Claim}

\theoremstyle{definition}

\DeclareMathOperator{\id}{id}

\title{The homotopy type of the independence complex of graphs with no induced cycles of length divisible by $3$}

\author{Jinha Kim}
\address{Discrete Mathematics Group\\
Institute for Basic Science (IBS)\\
Daejeon\\
Republic of Korea}
\email{jinhakim@ibs.re.kr}
\thanks{This work was supported by the Institute for Basic Science (IBS-R029-C1).}
\date\today

\begin{document}
\maketitle
\begin{abstract}
We prove Engstr\"{o}m's conjecture that the independence complex of graphs with no induced cycle of length divisible by $3$ is either contractible or homotopy equivalent to a sphere.
Our result strengthens a result by Zhang and Wu,
verifying a conjecture of Kalai and Meshulam which states that the total Betti number of the independence complex of such a graph is at most $1$.
A weaker conjecture was proved earlier by Chudnovsky, Scott, Seymour, and Spirkl, who showed that in such a graph, the number of independent sets of even size minus the number of independent sets of odd size has values $0$, $1$, or $-1$.
\end{abstract}

\section{Introduction}
We assume all graphs are finite and contain no loops and no multiple edges.
A subgraph of a graph $G$ is an {\em induced subgraph} if it can be obtained from $G$ by deleting vertices and all edges incident with those vertices.
An {\em induced cycle} is an induced subgraph that is a cycle.
An {\em independent set} is a set of pairwise non-adjacent vertices.
The {\em independence complex} of a graph $G$ is the abstract simplicial complex $I(G)$ on the vertex set $V(G)$ whose faces are the independent sets of $G$.
A graph is {\em ternary} if it contains no induced cycle of length divisible by $3$.

Here is our main theorem.
\begin{theorem}\label{main}
A graph is ternary if and only if every induced subgraph has the independence complex that is contractible or homotopy equivalent to a sphere.
\end{theorem}
We can easily deduce the converse of Theorem~\ref{main} as follows.
Kozlov \cite[Proposition 5.2]{Koz99} showed that the independence complex of a cycle is not homotopy equivalent a sphere if and only if the cycle has length divisible by $3$.
More precisely, if $C_\ell$ is a cycle of length $\ell \geq 3$, then the homotopy type of the independence complex is given by
\begin{align}\label{Koz}
    I(C_\ell) \simeq \begin{cases} S^{k} \vee S^{k} & \text{if }\ell=3k+3,\\S^{k} & \text{if }\ell=3k+2 \text{ or }3k+4.\end{cases}
\end{align}
Therefore, if every induced subgraph of a graph $G$ has the independence complex that is contractible or homotopy equivalent to a sphere, then $G$ does not contain an induced cycle of length divisible by $3$.

Our main result is motivated from a conjecture by Kalai and Meshulam.
For a simplicial complex $K$, let $\tilde{H}_i(K)$ be the {\em $i$-th reduced homology group} of $K$ over $\mathbb{Z}$ and $\tilde{\beta}_i(K)$ the {\em $i$-th reduced Betti number} of $K$, which is the rank of $\tilde{H}_i(K)$.
Let $\beta(K)$ be the {\em total Betti number}, that is, $\beta(K) = \sum_{i \geq 0}\tilde{\beta}_i(K)$.
Decades ago, Kalai and Meshulam~\cite{KM} conjectured that the independence complex of every ternary graph has total Betti number at most $1$.
This conjecture was recently proved by Zhang and Wu~\cite{WZ20}.

In general, for a graph $G$, $\beta(I(G)) \leq 1$ does not imply that $I(G)$ is contractible or homotopy equivalent to a sphere.
To see why, observe that the barycentric subdivision of any cell complex can be expressed as the independence complex of some graph.
For example, considering the real projective plane $\mathbb{R}\mathbb{P}^2$, one can find a graph whose independence complex is homotopy equivalent to $\mathbb{R}\mathbb{P}^2$, which has the total Betti number $0$, but is neither contractible nor homotopy equivalent to a sphere.

Our result, as well as the result by Zhang and Wu, is a generalization of a result about the reduced Euler characteristic of the independence complex of ternary graphs.
Given a simplicial complex $K$, the {\em reduced Euler characteristic} of $K$ is defined as \[\chi(K) = \sum_{i \geq 0} (-1)^i \tilde{\beta}_i(K).\]
It is a well-known fact in algebraic topology that $\chi(K) = \sum_{A \in K}(-1)^{|A|}$ (see \cite{Hat}).
Therefore, for a graph $G$, $|\chi(I(G))|$ is the difference between the number of independent sets of $G$ of even size and the number of those of odd size.
Kalai and Meshulam~\cite{KM} also posed a weaker conjecture that a graph $G$ is ternary if and only if $|\chi(I(H))| \leq 1$ for every induced subgraph $H$, and this conjecture was proved by Chudnovsky, Scott, Seymour and Spirkl~\cite{CSSS20}. 
\begin{theorem}\cite{CSSS20}\label{thm-char}
A graph $G$ is ternary if and only if $|\chi(I(H))| \leq 1$ for every induced subgraph $H$.
\end{theorem}

Earlier, Gauthier \cite{Gauthier} proved a special case of the conjecture: if a graph $G$ contains no (not necessarily induced) cycles of length divisible by $3$, then $|\chi(I(G))| \leq 1$.
By extending the work of Gauthier, Engstr\"{o}m~\cite{Eng20} showed that the independence complex of such a graph is either contractible or homotopy equivalent to a sphere.
Engstr\"{o}m conjectured that his result can be extended to ternary graphs.
Theorem~\ref{main} confirms Engstr\"{o}m's conjecture.

\medskip

Here is an overview of the paper.
In Section~\ref{sec:prelim}, we discuss some topological background on the homotopy type of simplicial complexes, including useful lemmas about the independence complexes of graphs. The proof of Theorem~\ref{main} will be presented in Section~\ref{sec:main}.

\section{Preliminaries}\label{sec:prelim}
In this section, we introduce some topological background (for details, see \cite{Hat}) and useful lemmas to determine the homotopy type of independence complexes.

For a graph $G$ and $v \in V(G)$, let $N(v) := \{u \in V: uv\in E\}$ and $N[v] := N(v) \cup \{v\}$.
For $W \subset V(G)$, let $N[W]:=\cup_{w \in W} N[w]$.

\subsection{Mayer-Vietoris Sequences}\label{sec-mvseq}
Let $K,A$ and $B$ be simplicial complexes such that $K=A\cup B$.
The {\em Mayer-Vietoris sequence} of the triple $(K,A,B)$ is the following long exact sequence of the homology groups of $K,A,B$ and $A \cap B$:
\[\cdots \to \tilde{H}_i (A \cap B) \to \tilde{H}_i(A) \oplus \tilde{H}_i(B) \to \tilde{H}_i(K) \to \tilde{H}_{i-1} (A \cap B) \to \cdots .\]

The following are basic observations about exact sequences.
\begin{itemize}
    \item If $0 \to A \to B \to 0$ is exact, then $A \simeq B$.
    \item If $0 \to A \to B \to C \to 0$ is exact, then $\text{rk}(B)=\text{rk}(A)+\text{rk}(C)$ where $A,B,C$ are finitely generated $\mathbb{Z}$-modules.\footnote{This was left as an exercise in p.146 of \cite{Hat}. The proof can be found, for example, in \cite{exact-rank}.}
\end{itemize}

Now let $G$ be a graph on $V$.
For each $v \in V$, observe that every independent set of $G$ containing $v$ is contained in $G - N(v)$.
This implies 
\[
I(G)=I(G-v) \cup I(G-N(v)).
\]
Note that $I(G-N(v))$ is a cone with apex $v$, thus it is contractible.
Finally, we observe
\[
I(G-v) \cap I(G-N(v))=I(G-N[v]).
\]
Then, by applying the Mayer-Vietoris sequence to the triple $(I(G), I(G-v), I(G-N(v)))$, we obtain the following long exact sequence:
\begin{align}\label{mv}
\cdots \to \tilde{H}_i (I(G-N[v])) \to \tilde{H}_i(I(G-v)) \to \tilde{H}_i(I(G)) \to \tilde{H}_{i-1} (I(G-N[v])) \to \cdots .    
\end{align}

For disjoint subsets $X$ and $Y$ of $V$ such that $X$ is independent in $G$, let $G(X|Y)$ be the subgraph of $G$ induced by $V-N[X]-Y$.
If $v \not\in N[X] \cup Y$, then we obtain the following exact sequence by replacing $G$ with $G-N[X]-Y$ in \eqref{mv}:
\begin{align}\label{mvseq}
\begin{split}
& \cdots \to \tilde{H}_{i}(I(G(X\cup\{v\}|Y))) \to \tilde{H}_i(I(G(X|Y\cup\{v\}))) \to \tilde{H}_i(I(G(X|Y))) \\
& \to  \tilde{H}_{i-1}(I(G(X\cup\{v\}|Y))) \to \tilde{H}_{i-1}(I(G(X|Y\cup\{v\}))) \to \tilde{H}_{i-1}(I(G(X|Y))) \to \cdots.
\end{split}
\end{align}

Recalling that $\tilde{H}_i(S^k)=0$ if $i \neq k$ and $\tilde{H}_k(S^k)\simeq \mathbb{Z}$, we can prove the following lemma.

\begin{lemma}\label{mvseq-betti}
Let $A,B$ and $C$ be simplicial complexes such that the following sequence is exact:
\[\cdots \to \tilde{H}_{i}(A) \to \tilde{H}_i(B) \to \tilde{H}_i(C) \to  \tilde{H}_{i-1}(A) \to \cdots.\]
Suppose $A \simeq S^k$ and $B \simeq S^\ell$ for some non-negative integers $k$ and $\ell$.
Then the following hold.
\begin{enumerate}[(i)]
    \item If $k>\ell$, then
$\tilde{\beta}_{k+1}(C)=\tilde{\beta}_\ell(C)=1$.
    \item If $k=\ell$, then either $\tilde{\beta}_i(C)=0$ for all non-negative integer $i$ or both $\tilde{H}_{k+1}(C) \neq 0$ and $\tilde{H}_{k}(C) \neq 0$.
\end{enumerate}
\end{lemma}
\begin{proof}
To prove (i), assume $k>\ell$.
Since $\tilde{H}_{\ell}(A)=\tilde{H}_{\ell-1}(A)=0$, we have $\tilde{H}_\ell(C) \simeq \tilde{H}_\ell(B)$.
Similarly, since $\tilde{H}_{k+1}(B)=\tilde{H}_{k}(B)=0$, we obtain
$\tilde{H}_{k+1}(C) \simeq \tilde{H}_{k}(A)$.
Therefore, $\tilde{\beta}_{k+1}(C)=\tilde{\beta}_k(A)=1$ and $\tilde{\beta}_\ell(C)=\tilde{\beta}_\ell(B)=1$.

Now to prove (ii), suppose $k=\ell$.
Since $\tilde{H}_{i-1}(A)=\tilde{H}_{i}(B)=0$ for $i \neq k,k+1$, we have $\tilde{H}_i(C)=0$ for $i \neq k,k+1$.
It is sufficient to show that if one of $\tilde{H}_k(C)$ and $\tilde{H}_{k+1}(C)$ is the trivial group, then the other has rank $0$.
Suppose $\tilde{H}_{k+1}(C)=0$.
Then we have a short exact sequence
\[0 \to \tilde{H}_{k}(A) \to \tilde{H}_k(B) \to \tilde{H}_k(C) \to 0.\]
This implies $\tilde{\beta}_k(B)=\tilde{\beta}_k(A)+\tilde{\beta}_k(C)$.
Since $\tilde{\beta}_k(A)=\tilde{\beta}_k(B)=1$ by the assumption, we have $\tilde{\beta}_k(C)=0$.
Applying a similar argument, one can show that $\tilde{\beta}_{k+1}(C)=0$ if $\tilde{H}_k(C)=0$.
\end{proof}

\subsection{Homotopy type theory}
Let $A$ and $B$ be two topological spaces.
\begin{itemize}
\item $A$ and $B$ are {\em homotopy equivalent} if there are continuous maps $f:A \to B$ and $g:B \to A$ such that $g \circ f \simeq \id_A$ and $f \circ g \simeq \id_B$, where $\id_X$ is the identity map on $X$.
We write $A \simeq B$ if $A$ and $B$ are homotopy equivalent.
In particular, if $A$ is contractible, i.e. $A$ is homotopy equivalent to a point, then we write $A \simeq *$.

\item The {\em wedge sum} of $A$ and $B$ is the space obtained by taking the disjoint union of $A$ and $B$ and identifying a point of $A$ and a point of $B$. We denote the wedge sum of $A$ and $B$ by $A \vee B$.

\item Let $\sim$ be an equivalence relation on $A$. Then we denote the quotient space of $A$ under $\sim$ by $A / \sim$.
Let $B \subset A$. Then we define $A / B$ as the quotient space $A / \sim$ where for all $a \neq b$ in $A$, $a \sim b$ if and only if $a, b \in B$.

\item The {\it suspension} of $A$ is the quotient space \[\Sigma A := A \times [0,1] / \sim\] where for all $(a, s) \neq (b, t)$ in $A \times [0,1]$, $(a, s) \sim (b, t)$ if and only if either $s=t = 0$ or $s=t = 1$.
\end{itemize}
Note that if $S^n$ is the $n$-dimensional sphere, then $\Sigma S^n \simeq S^{n+1}$.

Now let $K$, $K_1$ and $K_2$ be simplicial complexes where $K_1 \cap K_2 \neq \emptyset$, and let $L \neq \emptyset$ be a subcomplex of $K$.
Then,
\begin{enumerate}[(A)]
\item\label{aaaaaa} If $K=K_1 \cup K_2$, then $K/K_2 \simeq K_1/(K_1 \cap K_2)$.
\item\label{bbbbbb} Suppose the inclusion map $L \hookrightarrow K$ is homotopic to a constant map $c:L \to K$, that is, $L$ is contractible in $K$. Then $K/L \simeq K \vee \Sigma L$.
In particular, $K/L \simeq \Sigma L$ when $K$ is contractible, and $K/L \simeq K$ when $L$ is contractible.
\end{enumerate}

By applying \eqref{bbbbbb}, we can deduce the following well-known statement:
\begin{lemma}\label{homoto-quo}
Let $X$ be a simplicial complex, and $Y$ be a subcomplex of $X$.
If $X \simeq S^k$ and $Y \simeq S^\ell$ for some non-negative integers $k$ and $\ell$ with $\ell<k$, then $X/Y$ is homotopic to $S^k \vee S^{\ell+1}$.
\end{lemma}

By a similar argument as in Section~\ref{sec-mvseq}, we obtain the following lemma about homotopy equivalence of independence complexes.
\begin{lemma}\label{mv-homotopy-gen}
Let $G$ be a graph and $v$ a vertex of $G$.
If $X$ and $Y$ are disjoint subsets of $V(G)$ such that $X$ is independent, $N[X] \cup Y$ does not contain $v$, and $N[X] \cup N[v] \cup Y \neq V(G)$, then
\[
I(G(X|Y)) \simeq I(G(X|Y\cup\{v\}))/I(G(X\cup\{v\}|Y)).
\]
\end{lemma}

\begin{proof}
Observe that \[I(G(X|Y))=I(G(X|Y\cup\{v\}) \cup I(G(X|Y)-N(v))\] and \[I(G(X|Y\cup\{v\}) \cap I(G(X|Y)-N(v))=I(G(X\cup\{v\}|Y)).\]
Applying \eqref{aaaaaa}, we obtain \[I(G(X|Y))/ I(G(X|Y)-N(v)) \simeq I(G(X|Y\cup\{v\}))/I(G(X\cup\{v\}|Y)).\]
Since $I(G(X|Y)-N(v)) \simeq *$, by \eqref{bbbbbb}, we have \[I(G(X|Y)/ I(G(X|Y)-N(v)) \simeq I(G(X|Y)).\]
Therefore,
\[I(G(X|Y)) \simeq I(G(X|Y\cup\{v\}))/I(G(X\cup\{v\}|Y)),\]
as required.
\end{proof}

\section{Proof of Theorem~\ref{main}}\label{sec:main}
In this section, we prove the main result.
By \eqref{Koz}, it is sufficient to show the following.
	\begin{theorem}\label{thm:main2}
    Let $G$ be a ternary graph.
    Then $I(G)$ is either contractible or homotopy equivalent to a sphere.
	\end{theorem}

To prove Theorem~\ref{thm:main2} by contradiction, take a counter-example $G$ on $V$ which is minimal in the following sense: $I(G)$ is neither contractible nor homotopy equivalent to a sphere, but $I(H)$ is either contractible or homotopy equivalent to a sphere for every proper induced subgraph $H$ of $G$.

Let $X$ and $Y$ be disjoint vertex subsets of $G$.
We define $d(X|Y)$ as the following:
\[d(X|Y) = \begin{cases} d & \text{if }X \text{ is independent and }I(G(X|Y)) \simeq S^d,\\ * &\text{otherwise}.\end{cases}\]
If a graph $G$ has no vertex, then we write $I(G) \simeq S^{-1}$.
Note that if $X$ is an independent set, then $d(X|Y)=*$ implies that $I(G(X|Y))$ is contractible.

We need two lemmas.
The first lemma\footnote{Note that if $\tilde{\beta}_d(I(G(X|Y))) \neq 0$ for some $d$, then $I(G(X|Y)) \simeq S^d$ by the minimality of $G$, and hence $d(X|Y)=d$. In this sense, Lemma~\ref{triple} is a natural analogue of \cite[Lemma 3.1]{WZ20}. We include a proof for the completeness of the paper.} describes all possible types of the triples of the form $(d(X|Y),d(X\cup\{v\}|Y),d(X|Y\cup\{v\}))$ under certain conditions, and the second lemma shows that there is a universal constant $k$ such that $d(\emptyset|v) = d(v|\emptyset) = k$ for all $v$.
\begin{lemma}\label{triple}
Let $X$ and $Y$ be vertex subsets of $G$ such that $X \cup Y \neq \emptyset$ and $X \cap Y =\emptyset$.
For every vertex $v \not\in X \cup Y$, the triple $(d(X|Y),d(X\cup\{v\}|Y),d(X|Y\cup\{v\}))$ equals to one of the following:
\[(*,*,*), (k,*,k), (*,k,k), (k+1,k,*)\]
for some integer $k \geq -1$.
\end{lemma}
\begin{proof}
If $X\cup Y = V$, then there is nothing to prove, so we assume $X \cup Y \neq V$.
If $X$ is not an independent set, then we have $(d(X|Y),d(X\cup\{v\}|Y),d(X|Y\cup\{v\}))=(*,*,*)$.
Thus we assume that $X$ is an independent set.

If $v \in N(X)$, then we have $d(X\cup\{v\}|Y)=*$ since $X\cup\{v\}$ is not independent, and we have $G(X|Y) = G(X|Y\cup\{v\}))$.
Thus it must be $d(X\cup\{v\}|Y)=*$ and $d(X|Y)=d(X|Y\cup\{v\})$, implying that $(d(X|Y),d(X\cup\{v\}|Y),d(X|Y\cup\{v\}))$ is either $(*,*,*)$ or $(k,*,k)$ for some integer $k \geq -1$.

\medskip

Now assume $v \notin N(X)$.
There are two cases.

\smallskip

{\bf Case 1.} $N[X]\cup N[v] \cup Y=V$.
In this case, $G(X\cup\{v\}|Y)$ has no vertex, and hence $d(X\cup\{v\}|Y)=-1$.

If $N[X]\cup\{v\}\cup Y = V$, then $G(X|Y\cup\{v\})$ has no vertex and $G(X|Y)$ is a graph on $\{v\}$.
Since $I(G(X|Y))$ is a simplex on $\{v\}$, we have $(d(X|Y),d(X\cup\{v\}|Y),d(X|Y\cup\{v\}))=(*,-1,-1)$.

Thus, we may assume $N[X] \cup \{v\} \cup Y \neq V$.
Then we use the observation that the vertex set of $G(X|Y\cup\{v\})$ is contained in $N(v)$.
Since $G$ is ternary, $N(v)$ is an independent set of $G$.
This implies that $I(G(X|Y\cup\{v\}))$ is contractible, and hence $d(X|Y\cup\{v\})=*$.
Furthermore, we know that the vertex set of $G(X|Y)$, say $W$, is contained in $N[v]$.
Let $W' = W \setminus \{v\} \neq \emptyset$.
Then $G(X|Y)$ is the complete bipartite graph on $\{v\} \cup W'$, whose independence complex is the disjoint union of a simplex on $\{v\}$ and a simplex on $W'$.
Therefore, we have $I(G(X|Y)) \simeq S^0$, so $(d(X|Y),d(X\cup\{v\}|Y),d(X|Y\cup\{v\})) = (0,-1,*)$.

\smallskip

{\bf Case 2.} $N[X]\cup N[v] \cup Y \neq V$.
Note that, by the assumption, we have $d(X\cup\{v\}|Y) \neq -1$ and $d(X|Y\cup\{v\}) \neq -1$.
By Lemma~\ref{mv-homotopy-gen}, we obtain
\[I(G(X|Y)) \simeq I(G(X|Y\cup\{v\}))/I(G(X\cup\{v\}|Y)).\]

If $d(X\cup\{v\}|Y)=*$, then by the assumption, it is obvious that $I(G(X \cup \{v\}|Y))$ is contractible.
By \eqref{bbbbbb}, we know 
\[I(G(X|Y)) \simeq I(G(X|Y\cup\{v\})),\]
so $(d(X|Y),d(X\cup\{v\}|Y),d(X|Y\cup\{v\}))$ should be either $(*,*,*)$ or $(k,*,k)$ for some non-negative integer $k$.
Thus we may assume that $d(X\cup\{v\}|Y)=k$ for some non-negative integer $k$. 

If $d(X|Y\cup\{v\})=*$, then $I(G(X|Y\cup \{v\}))$ is contractible, and hence by \eqref{bbbbbb}, we have
\[I(G(X|Y)) \simeq \Sigma I(G(X\cup\{v\}|Y)).\] 
Since $I(G(X\cup\{v\}|Y)) \simeq S^k$, we obtain $I(G(X|Y)) \simeq S^{k+1}$, which implies $(d(X|Y),d(X\cup\{v\}|Y),d(X|Y\cup\{v\})) = (k+1,k,*)$.

Now suppose $d(X|Y\cup\{v\})=\ell$ for some non-negative integer $\ell$.
If $k<\ell$, then $I(G(X|Y)) \simeq S^\ell \vee S^{k+1}$ by Lemma~\ref{homoto-quo}, which is a contradiction to the minimality assumption of $G$.
Thus we may assume $k \geq \ell$.
If $k>\ell$, then $\tilde{\beta}_{k+1}(I(G(X|Y))=\tilde{\beta}_\ell(I(G(X|Y))=1$ by Lemma~\ref{mvseq-betti} (i), which is again a contradiction to the minimality of $G$.
Therefore, it must be $k=\ell$.
Now by Lemma~\ref{mvseq-betti} (ii), $I(G(X|Y))$ either has the zero Betti number in all dimensions or has non-trivial homology groups in two different dimensions.
By the assumption, we know that $I(G(X|Y))$ is either contractible or homotopy equivalent to a sphere.
Combining the above information, we conclude that $I(G(X|Y))$ should be contractible, and $(d(X|Y),d(X\cup\{v\}|Y),d(X|Y\cup\{v\})) = (*,k,k)$.
\end{proof}

For the second lemma, we first need to have $\beta(I(G)) \leq 1$.
This can be obtained from Theorem~\ref{thm-char} and the following claim.
\begin{claim}\label{claim-betti}
If $\beta(I(G)) \geq 2$, then there is an integer $k \geq 0$ such that $\tilde{\beta}_k(I(G)) \geq 2$ and $\tilde{\beta}_i(I(G))=0$ for all $i \neq k$.
\end{claim}
Note that Claim~\ref{claim-betti} is a weaker version of \cite[Claim 3.3]{WZ20}. The proof can be found in Appendix~\ref{sec:appendix}.

If $\beta(I(G)) \geq 2$, then we have $|\chi(I(G))| \geq 2$ by Claim~\ref{claim-betti}.
Then, Theorem~\ref{thm-char} implies that $G$ is not a ternary graph, which is a contradiction.
Therefore, we have $\beta(I(G))~\leq~1$.

\begin{lemma}\label{claim-dim}
There is a non-negative integer $k$ such that $d(\emptyset|v)=d(v|\emptyset)=k$ for all $v \in V$.
\end{lemma}

\begin{proof}
Since $G$ is a ternary graph, if $N[v] = V$ for some $v \in V$, we have $I(G) \simeq S^0$, which is a contradiction to the assumption on $G$.
Thus, we may assume $N[v] \neq V$ for all $v \in V$.
Then by Lemma~\ref{mv-homotopy-gen}, for any vertex $v \in V$, we have
\[
I(G) \simeq I(G(\emptyset|v))/I(G(v|\emptyset)).
\]
Note that each of $I(G(\emptyset|v))$ and $I(G(v|\emptyset))$ is either contractible or homotopy equivalent to a sphere.
By \eqref{bbbbbb},
\begin{itemize}
    \item if $I(G(v|\emptyset)) \simeq *$, then $I(G) \simeq I(G(\emptyset|v))$, and
    \item if $I(G(\emptyset|v)) \simeq *$, then $I(G) \simeq \Sigma I(G(v|\emptyset))$.
\end{itemize}
In both cases, it is clear that $I(G)$ is either contractible or homotopy equivalent to a sphere.
Thus we may assume both $I(G(\emptyset|v))$ and $I(G(v|\emptyset))$ are homotopy equivalent to spheres.

Assume $d(v|\emptyset)=\ell$ and $d(\emptyset|v)=k$ for some non-negative integers $k$ and $l$.
If $k>\ell$, then $I(G) \simeq S^k \vee S^{\ell+1}$ by Lemma~\ref{homoto-quo}, which is a contradiction to $\beta(I(G)) \leq 1$.
If $k < \ell$, then by Lemma~\ref{mvseq-betti} (i), we have $\tilde{\beta}_k(I(G))=\tilde{\beta}_{\ell+1}(I(G))=1$, which implies $\beta(I(G)) \geq 2$. 
This is again a contradiction to $\beta(I(G)) \leq 1$.
Thus we conclude that $k=\ell$.

Now suppose there exist $u, v \in V$ such that $d(u|\emptyset)=d(\emptyset|u)=p$ and $d(v|\emptyset)=d(\emptyset|v)=q$ for two non-negative integers $p,q$ with $p<q$.
By Lemma~\ref{triple}, we have \[(d(v|\emptyset),d(u,v|\emptyset),d(v|u)) \text{ is either } (q,*,q) \text{ or } (q,q-1,*).\]
If $d(u,v|\emptyset)=q-1$, then \[(d(u|\emptyset),d(u,v|\emptyset),d(u|v))=(p,q-1,d(u|v))\] which is possible only when $p=q$ by Lemma~\ref{triple}.
Thus we obtain $$(d(v|\emptyset),d(u,v|\emptyset),d(v|u))=(q,*,q).$$
On the other hand, Lemma~\ref{triple} implies that it must be \[(d(u|\emptyset),d(u,v|\emptyset),d(u|v))=(p,*,p).\]
Combining the above information, we have \[(d(\emptyset|u),d(v|u),d(\emptyset|u,v))=(p,q,d(\emptyset|u,v)).\] However, Lemma~\ref{triple} implies $q = p-1$, which is a contradiction.
\end{proof}

Now we are ready to complete the proof of Theorem~\ref{thm:main2}.
By Lemma~\ref{claim-dim}, there is a non-negative integer $k$ such that $d(v|\emptyset)=d(\emptyset|v)=k$ for all $v \in V$.

We claim $d(u,v|\emptyset)=k-1$ for any two distinct vertices $u,v$ of $G$.
By Lemma~\ref{triple},
\[(d(v|\emptyset),d(u,v|\emptyset),d(v|u)) \text{ is either }(k,*,k) \text{ or } (k,k-1,*)\] and \[(d(\emptyset|u),d(v|u),d(\emptyset|u,v)) \text{ is either } (k,*,k) \text{ or } (k,k-1,*).\]
Then $d(v|u)$ must be $*$, and hence we obtain $$(d(v|\emptyset),d(u,v|\emptyset),d(v|u))=(k,k-1,*).$$
Since $d(u,v|\emptyset) = k-1 \neq *$, we obtain that $\{u,v\}$ is an independent set of $G$.
Since this holds for every pair of two distinct vertices $u$ and $v$, we conclude that the whole vertex set $V$ is an independent set of $G$.
Thus, $I(G)$ is contractible, which is a contradiction to the assumption on $G$.

Therefore, if $G$ is ternary, then $I(G)$ is either contractible or homotopy equivalent to a sphere.

\section*{Acknowledgement}
The author is grateful to Woong Kook for pointing out a mistake in the original proof.
The author thanks Sang-il Oum for his comments that significantly improved the presentation of the paper.
The author also thanks Marija Jeli\'{c} Milutinovi\'{c} for pointing out an error in the statement of Lemma~\ref{mvseq-betti} in the initial version of the paper.
The author also appreciate the anonymous reviewers for helpful comments.

\begin{appendix}
\section{Proof of Claim~\ref{claim-betti}}\label{sec:appendix}
As it was mentioned in Section~\ref{sec:main}, Claim~\ref{claim-betti} follows from \cite[Claim 3.3]{WZ20}.
Here, we give a proof of Claim~\ref{claim-betti} for the completeness of the paper.
We need the following lemma about an exact sequence of homology groups.
\begin{lemma}\label{lemma-appendix}
Let $A,B$ and $C$ be simplicial complexes such that the following sequence is exact:
\[\cdots \to \tilde{H}_{i}(A) \to \tilde{H}_i(B) \to \tilde{H}_i(C) \to  \tilde{H}_{i-1}(A) \to \cdots.\]
Suppose that each of $A$ and $B$ is either contractible or homotopy equivalent to a sphere.
If $\tilde{\beta}_k(C) >0$ and $\tilde{\beta}_{\ell}(C) >0$ for some distinct non-negative integers $k,\ell$, then either $A \simeq S^{k-1}$ and $B \simeq S^{\ell}$ or $A \simeq S^{\ell-1}$ and $B \simeq S^k$.
\end{lemma}
\begin{proof}
Since $\tilde{\beta}_k(C) \neq 0$, we know that either $\tilde{H}_{k-1}(A) \neq 0$ or $\tilde{H}_k(B) \neq 0$.
Similarly, since $\tilde{\beta}_{\ell}(C) \neq 0$, we have either $\tilde{H}_{\ell-1}(A) \neq 0$ or $\tilde{H}_{\ell}(B) \neq 0$.
Suppose $\tilde{H}_{k-1}(A) \neq 0$.
Then we have $A\simeq S^{k-1}$, and hence we obtain $\tilde{H}_{\ell-1}(A)=0$. This implies $\tilde{H}_{\ell}(B) \neq 0$.
Thus we can conclude that $A \simeq S^{k-1}$ and $B \simeq S^{\ell}$.
By a similar argument, we can obtain $A \simeq S^{\ell-1}$ and $B \simeq S^k$ if $\tilde{H}_k(B) \neq 0$.
\end{proof}

Now we are ready to prove Claim~\ref{claim-betti}.
\begingroup
\def\thetheorem{\ref{claim-betti}}
\begin{claim}
If $\beta(I(G)) \geq 2$, then there is an integer $k \geq 0$ such that $\tilde{\beta}_k(I(G)) \geq 2$ and $\tilde{\beta}_i(I(G))=0$ for all $i \neq k$.
\end{claim}
\addtocounter{theorem}{-1}
\endgroup
\begin{proof}
To show by a contradiction, assume that $\beta_k(I(G)) >0$ and $\beta_{\ell}(I(G)) >0$ for some non-negative integers $k$ and $\ell$ with $k < \ell$.
Take a vertex $v \in V(G)$.
By applying Lemma~\ref{lemma-appendix} to \eqref{mv}, we have either $d(v|\emptyset)=k-1$ and $d(\emptyset|v)=\ell$ or $d(v|\emptyset)=\ell-1$ and $d(\emptyset|v)=k$.
Thus we can partition the vertex set $V(G)$ into two parts $V_1$ and $V_2$, where
\[
V_1=\{v \in V(G): d(v|\emptyset)=k-1 \text{ and } d(\emptyset|v)=\ell\},
\]
\[
V_2=\{v \in V(G): d(v|\emptyset)=\ell-1 \text{ and } d(\emptyset|v)=k\}.
\]

First, we claim that $V_2$ is an independent set.
It is sufficient to show that any two vertices in $V_2$ are not adjacent.
Take $u,v \in V_2$.
By Lemma~\ref{triple}, $$(d(u|\emptyset),d(u,v|\emptyset),d(u|v)) \text{ is either }(\ell-1,\ell-2,*) \text{ or }(\ell-1,*,\ell-1)$$ and
$$(d(\emptyset|v),d(u|v),d(\emptyset|u,v)) \text{ is either }(k,k-1,*) \text{ or }(k,*,k).$$
If $d(u|v) \neq *$, then $k-1=d(u|v)=\ell-1$, which is a contradiction to $k<\ell$.
Thus we have $$(d(u|\emptyset),d(u,v|\emptyset),d(u|v))=(\ell-1,\ell-2,*).$$
Then $d(u,v|\emptyset)=\ell-2$ implies that $u$ and $v$ are not adjacent, as required.

Now, if $V_1 = \emptyset$, then $I(G)$ is a simplex on $V_2 \neq \emptyset$, which is contractible.
This is a contradiction to the assumption.
Therefore, we may assume $V_1 \neq \emptyset$.

Next, we claim that for any two disjoint subsets $X$ and $Y$ of $V_1$ such that $X \cup Y \neq \emptyset$,
\begin{align}\label{triple-induction}
    d(X|Y)= \begin{cases} k-|X| & \text{if }|Y|=0,\\$*$ & \text{if }|X|,|Y|>0,\\\ell &\text{if } |X|=0.\end{cases}
\end{align}
We prove by induction on $|X|+|Y|$.
If $|X|+|Y|=1$, then it is obvious from the definition of $V_1$.
Now suppose 
\begin{enumerate}[(i)]
\item \eqref{triple-induction} holds for any two disjoint subsets $X,Y$ of $V_1$ such that $|X|+|Y|=m$ for some positive integer $m < |V_1|$.
\end{enumerate}
Take $A \subseteq V_1$ with $|A|=m+1$.
By the induction hypothesis~(i), we have $d(\emptyset|A\setminus\{a\})=\ell$ for all $a \in A$.
Hence by Lemma~\ref{triple},  for all $a \in A$, we know that $$(d(\emptyset|A\setminus\{a\}),d(\{a\}|A\setminus\{a\}),d(\emptyset|A)) \text{ is either }(\ell,\ell-1,*) \text{ or } (\ell,*,\ell).$$

Now, for any partition $A=A_1 \cup A_2$ such that $|A_1|,|A_2|>0$, we claim that 
\begin{align}\label{induction-partition}
    d(A_1|A_2)=\begin{cases} \ell-1 & \text{if }d(\emptyset|A)=*,\\ $*$& \text{if }d(\emptyset|A)=\ell.\end{cases}
\end{align}
This can be shown by induction on $|A_1|$.
If $|A_1|=1$, then the statement is true by the above argument.
Suppose 
\begin{enumerate}[(ii)]
    \item \eqref{induction-partition} holds for any partition $A=A_1 \cup A_2$ with $|A_1|=n$ for some positive integer $n <m$.
\end{enumerate}
Take a partition $A=A_1\cup A_2$ with $|A_1|=n+1 \leq m$.
Then for any $a \in A_1$, since $|A_1\setminus\{a\}| = n > 0$ and $|A_2| = |A| - |A_1| = m-n > 0$, we have
$d(A_1\setminus\{a\}|A_2)=*$ by the induction hypothesis~(i).
By Lemma~\ref{triple}, this implies $d(A_1|A_2)=d(A_1\setminus\{a\}|A_2\cup\{a\})$.
By the induction hypothesis~(ii), we obtain
\begin{align*}
    d(A_1|A_2)=d(A_1\setminus\{a\}|A_2\cup\{a\})=\begin{cases} \ell-1 & \text{if }d(\emptyset|A)=*,\\ $*$& \text{if }d(\emptyset|A)=\ell.\end{cases}
\end{align*}
This shows that \eqref{induction-partition} holds for any partition $A=A_1 \cup A_2$ such that $|A_1|,|A_2|>0$.

Now take $a \in A$.
By the induction hypothesis~(i), we have $d(A\setminus\{a\}|\emptyset)=k-|A|+1$.
Then, by Lemma~\ref{triple}, $(d(A\setminus\{a\}|\emptyset),d(A|\emptyset),d(A\setminus\{a\}|\{a\}))$ equals to either $(k-|A|+1,k-|A|,*)$ or $(k-|A|+1,*,k-|A|+1)$.
Thus $d(A\setminus\{a\}|\{a\})$ must be either $*$ or $k-|A|+1$.

On the other hand, \eqref{induction-partition} implies that
\[d(A\setminus\{a\}|\{a\}) =\begin{cases} \ell-1 & \text{if }d(\emptyset|A)=*,\\ $*$& \text{if }d(\emptyset|A)=\ell.\end{cases}\]
If $d(A\setminus\{a\}|\{a\}) \neq *$, then $k-|A|+1=d(A\setminus\{a\}|\{a\})=\ell-1$, which is impossible since $k-|A|+1=k-m \leq k-1 <\ell-1$.
Thus we must have $d(A\setminus\{a\}|\{a\})=*$.
This implies $d(A|\emptyset)=k-|A|$ and $d(\emptyset|A)=\ell$.
In addition, by \eqref{induction-partition}, $d(\emptyset|A)=\ell$ implies $d(A_1|A_2) = *$ for any partition $A = A_1 \cup A_2$ with $|A_1|,|A_2| > 0$.
This shows that \eqref{triple-induction} holds for any two disjoint subsets $X$ and $Y$ of $V_1$ such that $X \cup Y \neq \emptyset$.

Recall that $V_2$ is an independent set.
This implies that $I(G(\emptyset|V_1))$ is a simplex on $V_2$. 
Thus,
\[d(\emptyset|V_1)=\begin{cases}*&\text{ if }V_2\neq\emptyset,\\-1&\text{ if }V_2=\emptyset.\end{cases}\]
On the other hand, we have $d(\emptyset|V_1)=\ell$ by \eqref{triple-induction}.
This is a contradiction, since $\ell$ is a non-negative integer.

Therefore, there exists a non-negative integer $k$ such that $\tilde{\beta}_k(I(G))=\beta(I(G))\geq2$ and $\tilde{\beta}_i(I(G))=0$ for all $i \neq k$.
This completes the proof.
\end{proof}
\end{appendix}


\begin{thebibliography}{alpha}
\bibitem{CSSS20}
M. Chudnovsky, A. Scott, P. Seymour and S. Spirkl.
\newblock Proof of the Kalai-Meshulam conjecture.
\newblock {\em Isr. J. Math.}, 238:639--661, 2020.

\bibitem{Eng20}
A. Engstr\"{o}m.
\newblock On the topological Kalai-Meshulam conjecture.
\newblock {\em arXiv:2009.11077}, 2020.

\bibitem{Gauthier}
G. Gauthier.
\newblock The structure of graphs with no cycles of length $0$ (mod $3$).
\newblock {\em PhD thesis}, Princeton University, 2017.

\bibitem{Hat}
A. Hatcher, {\em Algebraic Topology}, Cambridge University Press, 2009.

\bibitem{KM}
G. Kalai.
\newblock When do a few colors suffice?
\newblock \url{https://gilkalai.wordpress.com/2014/12/19/when-a-few-colors-suffice}.

\bibitem{Koz99}
D.~N.~Kozlov.
\newblock Complexes of directed trees.
\newblock {\em J. Combin. Theory Ser. A}, 88(1):112--122, 1999.

\bibitem{WZ20}
W. Zhang and H. Wu.
\newblock The Betti Number of the Independence Complex of Ternary Graphs.
\newblock {\em arXiv:2011.10939}, 2020.

\bibitem{exact-rank}
\newblock \url{https://math.stackexchange.com/questions/1641562/exact-sequences-of-modules-and-rank}.

\end{thebibliography}
\end{document}